\input AHTOHFIE.STY
%\baselineskip 11.5pt
%\tolerance10000
%\hfuzz11pt
%\emergencystretch1mm
%%%%%%%%%%%%%%%%%%%%%%%%%%%%%%%%%%%%%%%%%%%%%%%%%%%%%%%%%%%
%%%%%%%%%%%%%%%%%%%%%%%%%%%%%%%%%%%%%%%%%%%%%%%%%%%%%%%%%%%
%%%%%%%%%%%%%%%%%%%%%%%%%%%%%%%%%%%%%%%%%%%%%%%%%%%%%%%%%%%

\UDC{
512.542     %Конечные группы
+ 512.54.03 %Элементарные теории различных классов групп
}


\title{%
How many tuples of group elements have a given property?
}
\author{%
Anton A. Klyachko
\quad
Anna A. Mkrtchyan
}
\address{
\myAddress
\quad anna.mkr@gmail.com
}

{\center{\ss with an appendix by Dmitrii V. Trushin}
%\medskip
}

%{\vskip-\medskipamount\smallskip\ssqi
%\center{
%The Einstein Institute of Mathematics. The Hebrew University of Jerusalem
%\\
%trushindima@yandex.ru
%}
%\medskip

\grants{\RFBR11-01-00945}

\abstract{%
%\narrower%\narrower
Generalising Solomon's theorem, C.~Gordon and F.~Rodriguez-Villegas
have proven recently that, in any group, the number of solutions to a
system of coefficient-free equations is divisible by the order of this
group whenever the rank of the matrix composed of the exponent sums of
$i$th unknown in $j$th equation is less than the number of
unknowns. We generalise this result in two directions: first, we
consider equations with coefficients, and secondly, we consider not
only systems of equations but also any first-order formulae in
the group language (with constants). Our theorem implies some
amusing facts; for example, the number of group elements whose squares
lie in a given subgroup is divisible by the order this subgroup.
}

%%%%%%%%%%%%%%%%%%%%%%%%%%%%%%%%%
\s 0.
Introduction

\proclaim Solomon theorem \rm [Solo69].
In any group, the
number of solutions to a system of coefficient-free equations is
divisible by the order of this group if the number of equations is
less than the number of unknowns.

This topic was developed in different directions
(see, e.g.
[Stru95],
[AmV11],
[Isaa70],
and references therein),
but the simplest and most natural generalisation
of this theorem was
obtained quite recently.

\proclaim Gordon--Rodriguez-Villegas theorem \rm \[GRV12].
In any group, the
number of solutions to a system of coefficient-free equations is
divisible by the order of this group if the rank of the matrix composed
of exponent sums of $i$th unknown in $j$th equation is less than the
number of unknowns.

For example, the Solomon theorem
does not apply
to the system $x^2y^3[x,y]y^{-1}=1=(yx)^2$,
but nevertheless, the number of its solutions
in a group
is divisible by the order of the group
according to the
Gordon--Rodriguez-Villegas theorem,
because the rank of the corresponding matrix
$
\left({\scriptscriptstyle2\;2 \atop
\scriptscriptstyle2\;2}\right)
$
is one while there are two unknowns.

What if the rank of the matrix is much less than the number of
unknowns? Does this imply that
the number of solutions must be divisible by a higher power of the order
of the group? The answer is \emph{no}. And even a
weaker conjecture is false. In Section 2, we give
an example.

We generalise the Gordon--Rodriguez-Villegas theorem in other directions.
We study equations with coefficients, and not only
systems of equations but also any first-order formulae.
The main theorem implies a lot of non-obvious facts, e.g. that mentioned
in the abstract. Section 1 contains the main theorem.
In Section 2, we give several examples.
In Section 3, we prove the main theorem.
In particular, our argument proves
the Gordon--Rodriguez-Villegas theorem and is (in this case)
somewhat simpler than the original one (from our point of view), although
our proof
is based on the same ideas. In Section 4, we give a direct proof of the
amusing fact from the abstract.
We have not succeeded in finding this fact in literature,%
\fn{\rm
In 2017, we learned that this fact was proven in [Iwa82].
}
though it might be easily
obtained from a result of P. Hall (see Section 1)
generalising the well-known Frobenius
theorem
[Frob03] (see also [Hall59]) which says that
{\sl
the number of solutions to the equation
$x^n=g$ is divisible by~$\GCD(n, |C(g)|)$}.
The Frobenius theorem was generalised in various directions
(see, e.g.~[Hall36],
[Kula38],
[Sehg62],
[BrTh88],
[AsTa01],
and references therein).

The appendix written by D. V. Trushin contains a proof of a pure
logical proposition
which makes it possible to simplify somewhat the statement of a corollary
of the main theorem at the expense of preliminarily transformation of the
logic formula.

\smallskip
\noindent
{\bf Our notation}
is mainly standard. Note only that if
$k\in \Z$ and $x$ and $y$ are elements of a group, then $x^y$,
$x^{ky}$, and~$x^{-y}$ denote $y^{-1}xy$, $y^{-1}x^ky$ and
$y^{-1}x^{-1}y$, respectively. A commutator~$[x,y]$ is
$x^{-1}y^{-1}xy$. If $X$ is a subset of a group, then
$|X|$, $\gp X$, and
$C(X)$
denote the cardinality of $X$, the
subgroup generated by $X$, and
the centraliser of~$X$.
The letter~$\Z$ denotes the set of integers.

%10.04.14

The authors thank an anonymous referee, A. V. Vasilev, and I. M. Isaacs
for useful remarks.

%%10.04.14

%%%%%%%%%%%%%%%%%%%%%%%%%%%%%%%%%
\s 1.
Main theorem

%10.04.14

Consider the group language $L$ over a group
$G$. This language has two functional symbols: $\cdot$ and~$^{-1}$, and
also, for each element of $G$, there is a constant symbol $g$.
We do not assume that the group is finite (although it is in the majority
of interesting cases); the assertions about divisibility should be
understood in the sense of cardinal arithmetics: any infinite cardinal
is divisible by any smaller or equal nonzero cardinal (and zero is
divisible by any cardinal).

Consider an arbitrary first-order formula $\phi$
in the language $L$.
Each atomic subformula can be written in the form
$$
u=1,
$$
where the words $u\in G*F$ and $F$ is the free group generated by all
(free and bound) variables of $\phi$.
Thus, the words $u$
(possibly, different for different subformulae)
can contain free and bound variables and elements
of~$G$ (called
\emph{coefficients} of $\phi$).

Now, we define the \emph{digraph $\Gamma(\phi)$ of the formula $\phi$} as
follows.
The vertices of $\Gamma(\phi)$ are all bound variables of~$\phi$. Each
atomic subformula containing bound
variables has the form
$$
v_1(y_1)w_1(x_1\dots,x_n)\dots v_r(y_r)w_r(x_1\dots,x_n)=h,
$$
where $y_i$ are bound variables of the formula
$\phi$ (not necessarily different), $x_1\dots,x_n$ are all (different)
free variables of~$\phi$, the words $v_i(y_i)$ belong to the free product
$G*\gp{y_i}_\infty$ of~$G$ and the infinite cyclic group
generated by the letter~$y_i$,
the words $w_i(x_1\dots,x_n)$ belong to the free product
$G*F(x_1\dots,x_n)$ of~$G$ and the free group with basis
$x_1\dots,x_n$, and $h\in G$.
Let us connect the vertices $y_i$ and~$y_{i+1}$ (subscripts modulo $r$)
by a directed edge labelled by an integer
tuple~$(\alpha_1,\dots,\alpha_n)$, where $\alpha_j$ is the exponent sum of
$x_j$ in $w_i$; loops labelled by zero
tuples are excluded.
We apply this construction to
each atomic subformula containing bound variables.

For instance, if the formula $\phi(x_1,x_2)$ has the form%
\fn{%
We do not assume that the formula is always in
the prenex normal form (i.e. all quantifiers are moved outside),
as in this example.
}
$$
\forall y \exists z\; \(
                        \([ygy,x_1gx_2]x_1z^{-1}azx_2^{-3}x_1^3hx_2^7=1\)
                        \wedge
                        \neg\(z^{-1}bz(x_2x_1)^2=1\)
                        \vee \((x_1^2x_2^2)^{5}=1\)
                       \),
\eqno{(1)}
$$
where $g,h,a,b\in G$ are some fixed elements (not necessarily different),
then the graph $\Gamma(\phi)$ looks as follows:%

\goodbreak
%\vskip1cm plus 1cm minus5mm
\bigskip
\centerline{\input 1.PIC}
\nobreak%
%\vskip5mm%
\centerline{Fig. \lowercase{1}}%
%\vskip1cm plus 1cm minus 5mm%
\goodbreak
\bigskip
%\par
%\noindent

Now, in the graph $\Gamma(\phi)$, we choose cycles $c_1,c_2,\dots$
generating the first homology group (e.g. we can take generators of the
fundamental groups of all components), and compose the
\emph{matrix $A(\phi)$ of~$\phi$} as follows: for each generating
cycle~$c_i$ we write a row which is the sum of labels of the edges of
this cycle
%10.13
(the edge labels are summed with signs plus or minus depending on
the orientation),
%10.13
and then we add rows consisting of exponent sums of atomic
subformulae containing no bound variables.

This matrix depends on the choice of generating cycles but
the integer
linear hull of its rows is determined uniquely by the formula
$\phi$. For the example above, the matrix $A(\phi)$ is
$$
A(\phi)=
\pmatrix{
-1&-1\cr
5&5\cr
2&2\cr
10&10\cr
}
\qbox{(for an obvious choice of three generating cycles).}
\eqno{(2)}
$$

A bound variable $t$ is called \emph{isolating}
if it occurs in atomic subformulae only in
subwords of the form  $t^{-1}g_it$, where $g_i\in G$.
The corresponding coefficients $g_i$ are called
\emph{isolated}. More precisely, an element $g$ of $G$ is called
\emph{isolated} if it occurs in $\phi$ only in
subwords o the form $t_i^{-1}gt_i$, where all $t_i$ are isolating
variables. In the example under consideration, $z$
is the only isolating variable and $a$ and $b$ are isolated coefficients.

\proclaim{Main theorem}.
If the rank of the matrix $A(\phi)$ of a formula $\phi$ is less than
the number of free variables of this formula, then
the number of tuples of group elements
satisfying formula $\phi$ is divisible by the order of
the centraliser
of the set of
all non-isolated coefficients of~$\phi$.
In particular, this number is divisible by the order of
the group if all
non-isolated
coefficients equal $1$.

In the example above, the rank of the matrix is one and there are two free
variables. Therefore, the theorem asserts that the cardinality of the set
$$
\left\{(x_1,x_2)\in G^2\;;\;
\forall y \exists z\; \(
                        \([ygy,x_1gx_2]x_1z^{-1}azx_2^{-3}x_1^3hx_2^7=1\)
                        \wedge
                        \neg\(z^{-1}bz(x_2x_1)^2=1\)
                        \vee \((x_1^2x_2^2)^{5}=1\)
                       \)
\right\}
$$
is divisible by
$|C({g,h})|$ (even if $g=a$, we should consider the coefficient $g$
as non-isolated). In the following section, we
give more instructive examples.

%%10.04.14

\Corollary 1.
The number of solutions to a system of equations in a group is divisible
by the order of the
centraliser
of the set of coefficients if the
rank of the matrix of this system is less than the number of unknowns.

This corollary transforms into the Gordon--Rodriguez-Villegas theorem in
the case where all coefficients equal 1.

\proclaim{Gordon--Rodriguez-Villegas conjugation theorem}
\rm(\[GRV12], Corollary 3.5%
\fn{translated into a language convenient for our purposes.}%
).
\newline
Let $\{w_j(x_1,\dots,x_n)\}\subset F(x_1,\dots,x_n)$ be
a set of words (elements of a free group) such
that the rank of the matrix composed of exponent sums of
$x_i$ in $w_j$ is less than $n$. Then, for any group~$G$
and any elements $h_j\in G$, the number of tuples satisfying the formula
$$
\bigwedge_j \(\exists q_j\; w_j(x_1,\dots,x_n)=q_j^{-1}h_jq_j\),
$$
is divisible by the order of~$G$.

This statement
(generalizing a conjugation theorem from [Solo69])
is obviously stronger than
the Gordon--Rodriguez-Villegas theorem
from the introduction.
Our theorem gives even stronger fact: the pointwise conjugation
is replaced by a common conjugation.

\Corollary 2.
Under the condition of the Gordon--Rodriguez-Villegas conjugation theorem,
the number of tuples satisfying the formula
$$
\exists q\;\(\bigwedge_j \(w_j(x_1,\dots,x_n)=q^{-1}h_jq\)\),
$$
is divisible by the order of~$G$.
\newline
\(If $w_j(x_1,\dots,x_n)\in G*F(x_1,\dots,x_n)$, then
the number of satisfying tuples is divisible by the order of
the centraliser of the set of all coefficients of the words $w_j$.\)

%10.04.14

\Proof
The graph contains the only vertex corresponding to an isolated
variable $q$;
coefficients $h_j$ are isolated;
and the matrix of the formula coincides with the matrix from
the Gordon--Rodriguez-Villegas conjugation theorem. So, Corollary 2
follows from the main theorem.

\medskip

The next proposition generalises the
Gordon--Rodriguez-Villegas conjugation theorem
in another direction.

\proclaim{Invariant set theorem}.
Suppose that
$U_j,V_j$ are subsets of
a finite group $G$ that are conjugation invariant
(i.e. they are unions of some conjugation classes)
and
$\{w_j(x_1,\dots,x_n)\}\subset F(x_1,\dots,x_n)*G$ is an
arbitrary set of words such that
the rank o the matrix composed of exponent sums of
$x_i$ in $w_j$ is less than $n$.
Then the number of tuples of elements of $G$
satisfying the formula
$$
\bigwedge_j \(w_jU_j\subseteq V_j\),
$$
is divisible by the order of
the centraliser of the set of all coefficients of all words
$w_j$.
In particular this number is a multiple of
$|G|$ if
$\{w_j(x_1,\dots,x_n)\}\subset F(x_1,\dots,x_n)$.

\Proof
Consider one inclusion $wU\subseteq V$.
Let us decompose $U$ and $V$ into unions of conjugation classes:
$$
U=a_1^G\cup a_2^G\cup\dots,
\quad
V=b_1^G\cup b_2^G\cup\dots.
$$
The inclusion $wU\subseteq V$
is equivalent to the following first-order formula:
$$
\forall y_1 \forall y_2\dots \exists z_1 \exists z_2\dots
\bigwedge_k\bigvee_l wy_k^{-1}a_ky_k=z_l^{-1}b_lz_l.
$$
Thus, the conjunction of inclusions is transformed into
a first-order formula where all bound variables are isolating and the
coefficients
$a_i$ and $b_i$ are isolated.
The matrix of this formula differ from the exponent-sum matrix
only in repeated rows and the assertion follows from the main theorem.

\medskip

This theorem transforms into
the Gordon--Rodriguez-Villegas conjugation theorem
in the case where $U_j=\1$ and
$w_j\in F(x_1,\dots,x_n)$.
Note that the assertion of the invariant set theorem remains true
if we replace the conjunction of inclusions
by the conjunction of arbitrary (maybe different)
set-theoretic relations ``logically expressible"
via inclusions. For instance,
$\subseteq$ may be replaced by
$\subset$, $\supseteq$, $\supset$, $=$, $\ne$, $\not\subset$,
``$intersects$",\dots
The conjunction itself may be replaced by any
quantifier-free first-order formula
For example, {\sl if $A=a^G$ and $B=b^G$ are some conjugation classes
of a group
$G$, then the number of  pairs $(x,y)\in G^2$ such that
$$
x^2y^3[x,y]y^{-1}AB=A
\qqbox{or}
(yx)^2A \hbox{ does not intersects } B,
$$
is divisible by $|G|$.
}
This follows from the main theorem. (The second term of the disjunction
is equivalent to
the formula
$\forall z \forall t\ (yx)^2z^{-1}az\ne t^{-1}bt$.)

%%10.04.14

The following statement is an analogue of the Solomon theorem
(and transforms into it in the case when all coefficients equal 1,
formula is quantifier-free and is a conjunction of equalities).

%10.04.14

\Corollary 3.
The number of tuples of group elements
satisfying formula $\phi$ is divisible by the order of the
centraliser of the
set
of all non-isolated coefficients of~$\phi$
(in particular, this number is divisible by the order of
the group provided
all non-isolated coefficients equal {\rm 1}) if, in the formula $\phi$,
$$
\eqalign{
(the\ number\ of\ proper\ occurrences\ of\ bound\ variables)+&\cr
+(the\ number\ of\ components\ of\ \Gamma(\phi))+&\cr
+(the\ number\ of\ atomic\ subformulae\
containing\ no\ bound\ variables)&
<
(the\ number\ of\ all\ variables).
}
\eqno{(*)}
$$

%%10.04.14

An \emph{occurrence of a variable $y$} is a maximal
subword in the left-hand side of an equation (the left-hand side
is considered as a cyclic word) containing the variable $y$ and no
other variables; An occurrence is called \emph{proper} if it
does not coincide with the entire left-hand side of the equation.  In the
example above, there are two occurrences of~$y$ and two occurrences of~$z$,
that is four occurrences of bound variables; all these occurrences are
proper. The number of occurrences of the special variable $q$ is always
equal to the number of inhomogeneous equations (because we
consider left-hand sides of equations as cyclic words).

\Proof
The rank of the first homology group of a graph equals to the number
of
edges minus the number of vertices plus the number of connected
components. The number of vertices equals to the number of bound
variables, the number of edges equals to the number of occurrences of
bound variables; improper occurrences give loops with zero labels.
Therefore, the rank of the matrix $A(\phi)$ is an most the left-hand side
of $(*)$ minus the number of bound variables. It remains to apply the
theorem.

\medskip

It is interesting that actually the graph $\Gamma(\phi)$ can always be
assumed to be
connected, in the sense that the graph of a formula written economically,
i.e. with minimal possible number of bound variables (among formulae
equivalent to the given one) is always connected. For example,
any formula with two bound variables that do not occur together in
atomic subformulae is equivalent to a formula with one bound variable.
Dima Trushin proves this fact in the appendix.

\Corollary 4.
The number of group elements whose $k${-th} powers belong to a given
subgroup is divisible by the order of this subgroup.%
\fn{\rm
In 2017, we learned that this fact was proven in [Iwa82].
}

\Proof
Let $H$ be a subgroup of a group $G$. We are interested in elements $x$
such that
$x^k\in H$. First, suppose that the subgroup $H$ is the centraliser of
a subgroup $A$. Then the inclusion $x\in H$ is equivalent to the system of
equations $\{[x^k,a]=1\;;\;a\in A\}$ satisfying the main theorem
(there are no
bound variables and the matrix is zero). Therefore, the number of
solutions is divisible by the order of the centraliser of the set of
coefficients, i.e. by $|H|$, as required.

Now let $H$ be an arbitrary subgroup.
Let us use the following trick.
We embed the group $G$ in
a larger group $\^G$ in such a way that $H$ becomes the centraliser of
a
subgroup $A$ of $\^G$. In addition, we must guarantee that
all solutions of our system of equations over $\^G$ belong to $G$.
For this sake, we make $G$ to be the centraliser of another subgroup
$B\subset\^G$ and consider the system of equations
$$
\(\bigwedge_{a\in A}\([x^k,a]=1\)\)
\wedge
\(\bigwedge_{b\in B}\([x,b]=1\)\).
$$
This will prove Corollary 4 in the general case.

The r\^ole of $\^G$ can be played by the amalgamated free product
$\^G=(B\times G)\zvezda_H (H\times A)$, where $A$ and $B$ are arbitrary
nontrivial centreless groups.  Clearly, $C(A)=H$ and
$C(B)=G$, i.e. the solutions to our system of equations are precisely
the elements of $G$ whose $k$th powers lie in~$H$. According to
the main theorem, the number of solutions is divisible by
$|C(A)\cap C(B)|=|H\cap G|=|H|$, as required.

This proof of Corollary 4 uses only a very particular case of the
main theorem, where
the formula $\phi$ is a system of equations with one
unknown.
This case of the main theorem follows immediately from an
old result of  P. Hall (generalising the Frobenius theorem).

\proclaim{Hall theorem \rm([Hall36], Theorem II)}.
In any group, the number of solutions to a system of equations
with one unknown is divisible by
$\GCD(|C|,n_1,n_2,\dots)$, where $C$
is the centraliser of the set of coefficients
and $n_i$ is the exponent sum of the unknown in
$i$-th equation.

A trick (another) allowing to
transform an arbitrary subgroup into a
centraliser can also be found in
[Hall36].
In Section 4, we give a direct proof of Corollary 4 illustrating a
part of the proof of the main theorem.

A more general fact
%10.13
%(which also might be derived from a known result [Sehg62])
can be obtained similarly.

\Corollary 5.
Suppose that $H$ is a subgroup of a group $G$ and $W$ is a subgroup (or a
subset) of a finitely generated
%10.13
%free
group $F$
%10.13
with infinite abelianisation $F/F'$.
Then the number of
homomorphisms $f\:F\to G$ such that $f(W)\subseteq H$ is divisible
by~$|H|$.

\Proof
%10.13
Consider a presentation
$F=\pres<X|R>$ of the group $F$.
The number of homomorphisms $f\:F\to G$ such that $f(W)\subseteq H$
coincides with
the number of solutions to the system of equations
$$
\(\bigwedge_{r\in R}\(r=1\)\)
\wedge
\(\bigwedge_{a\in A,\;w\in W}\([w,a]=1\)\)
\wedge
\(\bigwedge_{b\in B\;x\in X}\([x,b]=1\)\)
\qbox{with unknowns $X$}
$$
in the group
$\^G=(B\times G)\zvezda_H (H\times A)$,
where $A$ and $B$ are arbitrary
nontrivial centreless groups.
%10.13
The rank of the matrix of this system  coincides with the rank of
the matrix of the system $\{r=1\;;\; r\in R\}$
(as the other equations are commutators)
and is less than the number of unknowns $X$, because the
abelianisation of $F=\pres<X|R>$ is infinite.
According to
Corollary 1, the number of solutions is divisible by
$|C(A)\cap C(B)|=|H\cap G|=|H|$, as required.

\medskip

%10.13
Note that Corollary 5 coincides with Corollary 4 in the case where
$F=\Z$ and turns into the Gordon--Rodriguez-Villegas theorem when $H=G$.

%Следствие 4 можно слегка обобщить

%\Следствие 5.
%Число элементов группы, $k$-е степени которых лежат в данном смежном
%классе (левом или правом) по данной подгруппе, делится на порядок этой
%подгруппы.

%\Д
%Можно рассуждать аналогичным образом. Включение $x^k\in Hg$ равносильно
%системе уравнений $\{a^{x^k}=a^g\}$

% Нет! нельзя. В S_3 в (12)<(23)> квадрат только одного элемента лежит!

%%%%%%%%%%%%%%%%%%%%%%%%%%%%%%%%%
\s 2.
Examples

We start with a curious application of the Solomon theorem.

\Example 1.
We say that two elements of a group belong to the same \emph{tribe} if
their squares are equal. Clearly, the total size of all tribes is the
order of the group. It is less obvious that
\newline
\centerline{\sl
the sum of \the\year th powers of tribe sizes is a multiple of the
order of the
group.}
To prove this fact, it suffices to consider the system of equations
$x_1^2=\dots=x_{\the\year}^2$. Clearly, the number of solutions
is the sum of \the\year th powers of tribe sizes. The
number of equations is less than that of unknowns. So, the
statement is a corollary of the Solomon theorem.
The assertion remains
valid if \the\year\ is replaced by an arbitrary positive integer; the
squares (in the definition of tribes) can also be replaced by any (equal)
positive integer powers.

\Example 2.
{\sl
The number of pairs of group elements whose product of squares is a
cube is divisible by the order of the group.
}
This follows from
Corollary 3, because the formula $\exists z\; x^2y^2=z^3$ has one bound
variable, it occurs once, there are no equations without bound variables,
the graph is connected, and there are two free variables:
$1+0+1<2+1$.
\newline
This fact can also be derived from the
Gordon--Rodriguez-Villegas conjugation theorem.
Indeed, this theorem
implies that
the order of the group divides the number of pairs of group elements
whose product of squares is conjugate to any given element.
\newline
By the same reason,
{\sl the order of any group divides, e.g. the following
numbers:
\-
the number of pairs of noncommuting elements whose product of
squares is the cube of a noncentral element;
\-
the number of pairs of noncommuting
elements whose product of squares is a cube if and only if
the cube of their
product lies in the centre;
\-
the number of pairs of elements such that either
the product of their squares is a cube or their commutator
is not a square;
\-
\dots
}

\Example 3.
{\sl
The order of a group divides the
number of pairs of elements of this group whose product of squares is
the cube of a commutator \($x_1^2x_2^2=[z,t]^3$\) and
square of product is the commutator of
cubes of the same elements \($(x_1x_2)^2=[z^3,t^3]$\).
}
It is difficult to derive this fact from
the Gordon--Rodriguez-Villegas conjugation theorem, but it
follows immediately from our main theorem.
Corollary 2 gives
a stronger statement:
{\sl
the order of a group divides the number of pairs of
elements whose square of product and product of squares are
simultaneously conjugate to any given pair of elements.
}

The analogy between the Gordon--Rodriguez-Villegas theorem and well-known
properties of solutions to systems of linear equations
(over finite fields) could suggest an idea that, if the rank of
the matrix is much less than the number of unknowns, then the number
of solutions must be divisible by a higher power of the order of the
group. A more realistic question is the following.

{\sl
Is it true that the number of homomorphisms
from a finitely generated group $H$ into a group $G$ is divisible
by $|G|^m$
if $H$ admits an epimorphism onto a free
group of rank $m$?
}

\noindent
The point is that the number of solutions to a system of coefficient-free
equations
$$
\{u(x_1,\dots,x_n)=v(x_1,\dots,x_n)=\dots=1\}
$$
equals to the number of homomorphisms from the group
$H=\pres<x_1,\dots,x_n|u(x_1,\dots,x_n)=v(x_1,\dots,x_n)=\dots=1>$
to the group $G$. The matrix of the system has rank at most $r$ if and
only if the group $H$ admits an epimorphism onto the free abelian group
of rank $n-r$. The existence of an epimorphism onto the absolutely free
group of the same rank is a much stronger property, but
nevertheless, the conjecture under consideration is false for $m>1$,
as the following example shows.

\Example 4.
The group $\pres<x,y,z|z=z^xz^y>$ has an epimorphism onto
the
free group of rank two (sending $z$ to 1), but
the
number of solutions to the equation $z=z^{x}z^{y}$ in
the symmetric group $S_3$ is not divisible by $|S_3|^2=36$.
Indeed, with $z=1$,
there are 36 solutions ($x$ and $y$ can be arbitrary).
With $z=(123)$, there are $3\cdot3=9$ solutions
($x$ and~$y$ are arbitrary transpositions).
With $z=(321)$, there are also 9 solutions.
If $z$ is a transposition, then there are no solutions
(by parity). Thus, the total amount of solutions is $36+2\cdot9$.

%%%%%%%%%%%%%%%%%%%%%%%%%%%%%%%%%
\s 3.
Proof of the main theorem

\Lemma 1.
Under the conditions of the theorem, the first column of the matrix
$A(\phi)$ vanishes after a suitable invertible change of free variables.
In particular, the exponent sum of $x_1$ (in new variables) in each
cycle of the graph $\Gamma(\phi)$ is zero.

\Proof
The rank of the matrix $A(\phi)$ is less than the number of its columns;
therefore, some
integer (invertible) elementary transformations
of
columns produces a matrix with zero first column.  Elementary
transformations of columns are induced by obvious changes
of variables, e.g. the change $x_i \to x_ix_j^k$ adds $i$th column
multiplied by $k$ to $j$th column.  Lemma 1 is proven.

\medskip

For example, to annihilate the first column of matrix (2) from Section 1,
it suffices to subtract the second column from the first one, i.e.
the change of variables $x_2\to x_2x_1^{-1}$
transform formula (1) into the formula
$$
\forall y \exists z\; \(
                        \([ygy,x_1gx_2x_1^{-1}]x_1z(x_2x_1^{-1})^{-3}x_1^3h(x_2x_1^{-1})^7=1\)
                        \wedge
                        \neg\(z^8x_2^2=1\)
                        \vee \(x_1^2(x_2x_1^{-1})^2)^{5}=1\)
                       \).
\eqno{(3)}
$$
The graph of this formula is shown in Figure 2

\goodbreak
%\vskip1cm plus 1cm minus5mm
\bigskip
\centerline{\input 2.PIC}
\nobreak%
%\vskip5mm%
\centerline{Fig. \lowercase{2}}%
%\vskip1cm plus 1cm minus 5mm%
\goodbreak
\bigskip
%\par
%\noindent

and the matrix is
$$
\pmatrix{
0&-1\cr
0&5\cr
0&2\cr
0&10\cr
}.
$$
In what follows, we assume that the first column of $A(\phi)$
is zero.

%10.04.14

For each variable $t$ (free or bound), we
introduce new symbols $t^{(i)}$, where $i\in\Z$,
(their meaning will be
\newline
${t^{(i)}=t^{x_1^i}=x_1^{-i}tx_1^i}$
if $t$ is not isolating and
$t^{(i)}=t{x_1^i}$,
if $t$ is isolating
but we act formally at this step).
Similarly, for each coefficient $g$, we introduce new symbols
$g^{(i)}$, where $i\in\Z$.
Now, we transform formula $\phi$ as follows:
\item {1)}
in all atomic subformulae, replace each symbol
$t$ by $t^{(0)}$, where $t$ is a variable different from
$x_1$ or a coefficient;
\item {2)}
replace each subword of the form $(t^{(i)})^kx_1^l$, where
$t$ is a non-isolating variable or a coefficient, by $x_1^l(t^{(i+l)})^k$;
replace each subword of the form  $(g^{t^{(i)}})x_1^l$, where
$t$ is an isolating variable and $g\in G$, by
$x_1^l(g^{t^{(i+l)}})$;
\item {3)}
repeat step 2) while it is possible.

\enditem
This shift of symbols $x_1$ produces a formula without
$x_1$ (the exponent sum of $x_1$
vanishes in
each cycle of the graph and in each atomic subformula
without bound variables and, hence,
the exponent sum of $x_1$
vanishes in all atomic subformulae).
For example, formula (3) after these transformations will have
the following atomic subformulae:
$$
\eqalign{
[y^{(0)}g^{(0)}y^{(0)},g^{(-1)}x_2^{(-1)}]a^{z^{(-1)}}
\(x_2^{(-4)}x_2^{(-3)}x_2^{(-2)}\)^{-1}h^{(-7)}
x_2^{(-7)}x_2^{(-6)}x_2^{(-5)}x_2^{(-4)}x_2^{(-3)}x_2^{(-2)}x_2^{(-1)}&=1,
\cr
b^{z^{(0)}}(x_2^{(0)})^2&=1,
\cr
(x_2^{(-2)}x_2^{(-1)})^5&=1.
}
\eqno{(4)}
$$
Now, we proceed with transformation of the formula.

\item {4)}
In each homogeneous equation $\alpha$, we replace
all symbols
$t^{(i)}$ by $t^{(i+j_\alpha)}$, where
the integers~$j_\alpha$ are chosen such that,
for each bound variable $t$,
the symbols $t^{(i)}$ will occur in the entire formula
with at most one value of the superscript $(i)$;
this is possible, because the
exponent sum of $x_1$ vanishes in
each cycle of the graph.

\smallskip
\item{}
In formula (4), it
suffices to decrease the superscripts by one in the second equation.
In the general case, we can act as follows.
In each connected component $K$ of $\Gamma$, we
choose a vertex (variable) $y_K$. In each homogeneous equation $\alpha$
containing bound variables, we choose one of such variables
$y_\alpha$ and connect each vertex $y_\alpha$ by a path $p_\alpha$ with
the vertex~$y_K$ such that $y_\alpha\in K$.
The sum
$s_\alpha$ of
the first coordinates
of labels of edges of the path $p_\alpha$ does not depend on
the choice of the path by the condition.
Put $j_\alpha=-i_\alpha-s_\alpha$, where $i_\alpha$
is the unique number such
that $y_\alpha^{(i_\alpha)}$ occurs in the
equation~$\alpha$.
The
sum of the first coordinates
of labels of edges vanishes in each cycle, hence,
$s_\alpha=s_\beta$ if $y_\alpha=y_\beta$, the value $j_\alpha$ does not
depend on the choice of
variables $y_\alpha$ in equation~$\alpha$,
and the
changes $t^{(i)} \to t^{(i+j_\alpha)}$ in each
equation~$\alpha$
produce a formula such that
each bound variable $t$ occurs in the entire
formula only with one superscript
$(i)$, where $i$ is the sum of the first coordinates
of labels of
edges of any path from $y_K$ to $t$.
\item {5)}
Now, we replace each quantifier $\forall y$ and $\exists y$
by $\forall y^{(p)}$ and $\exists y^{(p)}$, where $p\in\Z$ is the
unique number such that $y^{(p)}$ occurs in atomic subformulae.
\item {6)}
Finally, we add equalities defining new symbols
to the obtained formula $\^\phi$, i.e.
we replace $\^\phi$ by the infinite formula
$$
\phi'=
\^\phi\wedge
\underbrace{
\(\bigwedge_{g}\(g^{(0)}=g\)\)
\wedge
\(\bigwedge_{t,i}\(t^{(i)}=x_1^{-1}t^{(i-1)}x_1\)\)
}_\delta,
\eqno{(**)}
$$
where $t$ ranges over all
free variables and all coefficients of the
initial formula, $i$ ranges over all integers,
and $g$ ranges over all coefficients.
The symbols $g^{(i)}$, where $g\in G$, are considered as free variables
of the formula $\phi'$.

\smallskip

\enditem
The numbers of tuples satisfying the
obtained formula $\phi'$ and the initial formula $\phi$
are equal.
Indeed, the formula~${\phi'=\^\phi\wedge\delta}$ admits
as many satisfying tuples as the formula
$\=\phi=\^\phi\big|_{t^{(i)}=t^{x_1^i}}$ (i.e. the formula $\^\phi$ with
each symbol~$t^{(i)}$, where $t$ is a free variable
of the initial formula or a
coefficient, is substituted by the expression $x_1^{-i}tx_1^i$).
The
formula~$\=\phi$ is equivalent to $\phi$ (i.e. this formulae
are satisfied by the same tuples). Indeed,
$\=\phi$ differs from~$\phi$ in two aspects:
\item{a)}
in quantifiers of $\=\phi$, symbols $y^{(p)}$
occurs instead of $y$;
\item{b)}
in atomic subformulae of $\=\phi$,
each bound variable $y$ is replaced by
the
expressions $(y^{(p)}){x_1^{-p}}$
or $(y^{(p)})^{x_1^{-p}}$ (depending on whether
$y$ is
isolating),
where $p$ is the same for all occurrences of $y$.

\enditem
Formulae different only in these aspects are obviously equivalent:
e.g.
$(\forall y\, \alpha(y,z,\dots))\equiv
(\forall t\, \alpha(t^{z^2},z,\dots))$,
because for any $g\in G$, if $y$ ranges over the whole group, then
$y^g$ ranges over the whole group.

%%10.04.14

\bigskip

Thus, it suffices to show that the number of tuples satisfying
$\phi'$ (such tuples are called \emph{solutions} henceforth),
is divisible by $|C|$, where
the letter $C$ denotes the centraliser
of the set
of
coefficients
of~$\phi$ (or of~$\phi'$, equivalently).

Consider a solution
$X=\(\~x_1,\,\~x_i^{(j)},\, \~g^{(j)}\; i=2,\dots,n,\; j\in\Z\)$.
The tuple $\(\~x_i^{(j)},\, \~g^{(j)}\)$,
i.e. everything but $\~x_1$, is
called the \emph{tail} of the solution
$X$.
Let $B_X$ denote the centraliser of the
tail of~$X$. Note that $B_X\subseteq C$
(because of the equations $g^{(0)}=g$ in formula $(**)$).

We say that two solutions are \emph{similar} if their tails are
conjugate by an element of $C$. Clearly, this is an equivalence
relation.
The theorem is a corollary of the following proposition.

\Proposition.
Each class of similar solutions contains exactly
$|C|$ solutions.

Let us find the number of solutions similar to $X$. The number of all
possible tails of such solutions is $|C|/|B_X|$, because, on the set of
tails of solutions similar to $X$, the group $C$ acts by conjugation
(since a tuple conjugate to the tail of a solution $Y$ by an element
$c\in C$ is also the tail of a solution, e.g., of $Y^c$) and $B_X$ is the
stabiliser of the tail of $X$.

The number of solutions with the same tail as that of~$X$ equals to
$|B_X|$, since, if a tuple with the same tail and with the first
coordinate $\~x_1'$ is also a solution, then the quotient
$\~x_1'(\~x_1)^{-1}$ must commute with the tail because of the
equations~$\delta$ in formula $(**)$, i.e.
$\~x_1'\in B_X\~x_1$.
On the other hand, any
element
$\~x_1'\in B_X\~x_1$ gives a solution with the same tail (as that of~$X$),
because the variable $x_1$ occurs in~$\phi'$ only in subformula $\delta$.

If $X'$ is a solution similar to $X$, then the number of solutions with
the same tail as that of~$X'$ equals to $|B_{X'}|=|B_X|$ (if tails are
conjugate, then their centralisers are conjugate and have the same order).

Thus, we obtain that the number of solutions similar to $X$ equals
$(|C|/|B_X|)\cdot|B_X|=|C|$, that proves the proposition and the theorem.

%%%%%%%%%%%%%%%%%%%%%%%%%%%%%%%%%
\s 4. Roots of subgroups

In Section 1, the following
assertion was derived from the main theorem.

\Corollary 4.
The number of group elements whose $k$-th powers belong to a given
subgroup is divisible by the order of this subgroup.

Here, we give a direct proof that exemplify
the concluding part of the proof of the main theorem.
For simplicity, we assume that $k=2$.
Let $H$ be a subgroup of a group $G$.
We are interested in elements $x\in G$ such that $x^2\in H$; such elements
are called \emph{solutions} henceforth. The assertion is implied by
the following lemma.

\Lemma.
Each double coset $HxH$ contains either 0 or $|H|$
solutions.

\Proof
Let $x$ be a solution; its \emph{tail} is the
coset $Hx$.

The group $H$ acts (on the right) on the set of tails of solutions
from the double coset $HxH$ by the right multiplication:
$$
Hy\circ h\:=Hyh\; (=the\ tail\ of\ the\ solution\ y^h).
$$
The stabiliser of the tail $Hx$ is $B_x\:=H\cap H^x$:
$$
Hx=Hxh \iff h\in H^x.
$$
Therefore, {\sl all possible solutions lying in $HxH$ have
precisely $|H|/|B_x|$ different tails.}

How many solutions have the same tail as that of $x$?
$$
Hx=Hy \imp yx^{-1}\in H,
$$
but if $y$ is also a solution, then
$$
(Hx)x=Hx^2=H=Hy^2=(Hx)y,
\qqbox{i.e.}
yx^{-1}\in H^x.
$$
Thus, each solution $y$ with the same tail as that of $x$ lies in
$B_xx$. On the other hand, each element from this coset
is a solution:
$$
(bx)^2=bxbx=bb^{x^{-1}}x^2\in H,
\qbox{because $b$, $b^{x^{-1}}$, and $x^2$ lie in $H$.}
$$
So, {\sl the number of solutions with the same tail as that of $x$ is
$|B_x|$}.

Since $|B_x|=|B_y|$ if $x$ and $y$ lie in the same double coset
$HxH$ (because $B_x$ and $B_y$ are conjugate in this case),
$HxH$ contains
${|B_x|\cdot(|H|/|B_x|)=|H|}$ solutions, that completes the
proof.

\medskip

Recently, I. M. Isaacs [Isaa12] obtained a character-theoretic proof of
this corollary.

In 2017, we learned that this fact was proven in [Iwa82].

%%%%%%%%%%%%%%%%%%%%%%%%%%%%%%%%%%%%%%%%
%%%%%%%%%%%%%%%%%%%%%%%%%%%%%%%%%%%%%%%%%
\bigskip\goodbreak\bigskip

{\ssdbf\center{\uppercase{
Appendix. On the minimisation of the number
of bound variables in first-order formulae
}}}
\medskip

\centerline{\ss Dmitrii V. Trushin}
\medskip
{{%\vskip-\medskipamount\smallskip
\ssqi
\center{
Einstein Institute of Mathematics,
The Hebrew University of Jerusalem,
\\
Givat Ram, Jerusalem, 91904, Israel
\\
trushindima@yandex.ru
}}
\medskip

Let $\phi$ be a first-order formula (in some language)
with bound variables $y_1,\ldots,y_k$
and free variables $x_1,\ldots,x_m$.
Consider the following graph $\Delta(\phi)$ with vertices
$y_1,\ldots,y_k$: vertices~$y_i$ and~$y_j$ are connected by an
edge if there exists an atomic subformula in $\phi$ containing
the both variables.

The following assertion shows that
the formula $\phi$
with disconnected graph $\Delta(\phi)$
is equivalent to a formula with fewer bound variables.

\Proposition.
Any formula $\phi$ is equivalent to a formula $\phi'$ with connected graph
$\Delta(\phi')$ such that
$$
|\Delta(\phi')|\leqslant\max(|\Delta_1|,\ldots,|\Delta_s|),
\qbox{
where
$\Delta_1,\ldots,\Delta_s$ are connected components of
$\Delta(\phi)$.
}
$$

\Proof
Let $\Phi_i$ be the set of formulae
$\psi$ such that
\item{1)}
all variables of $\psi$ (free and bound)
lie in the set $\{x_1,\ldots,x_m\}\cup\{y_1,\ldots,y_k\}$;
\item{2)}
all bound variables of $\psi$ lie in the
set $\{y_1,\ldots,y_k\}$;
\item{3)}
if a variable $y_j$ occurs in $\psi$ (as a bound or free variable),
then $y_j\in \Delta_i$.

\enditem
Note that the classes $\Phi_i$ are closed
under logical operators and quantifications on $y_j$.

Let $\Lambda$ be the closure of the union $\bigcup\limits_i\Phi_i$
with respect to logical operators ($\vee$, $\wedge$, and $\neg$).
Clearly, each formula from~$\Lambda$ can be written in the form
$$
\psi = \bigvee_{i=1}^l\bigwedge_{j=1}^s r_{ij},
\qbox{where $r_{ij}\in\Phi_j$},
\eqno{(D)}
$$
and the in form
$$
\psi = \bigwedge_{i=1}^l\bigvee_{j=1}^s r_{ij},
\qbox{where $r_{ij}\in\Phi_j$},
\eqno{(C)}
$$
because the conjunction and disjunction are mutually distributive,
the classes $\Phi_i$ are closed with respect to logical operations,
$\neg (A\vee B)=(\neg A)\wedge(\neg B)$, and
$\neg (A\wedge B)=(\neg A)\vee(\neg B)$.

\Lemma.
The class $\Lambda$ is closed with respect to
quantifications on the variables $y_j$.

\Proof
Note that
$$
\forall y (\psi_1(y) \wedge \psi_2(y)) =
(\forall \psi_1(y)) \wedge (\forall y \psi_2(y))
\quad \mbox{ and }\quad
\forall y (\psi_1(y) \vee \psi_2) =
(\forall y \psi_1(y))\vee \psi_2,
$$
where, in the second equality, $y$ is not a free variable of the formula
$\psi_2$. For the existence quantifier we have similar equalities
$$
\exists y (\psi_1(y) \vee \psi_2(y)) =
(\exists \psi_1(y))\vee (\exists y \psi_2(y))
\quad \mbox{ and }\quad
\exists y(\psi_1(y) \wedge \psi_2) = (\exists y \psi_1(y)) \wedge
\psi_2,
$$
where, in the second equality, $y$ is not a free variable of the formula
$\psi_2$.

Suppose that a formula $\psi\in\Lambda$ is written in the form $(D)$:
$$
\psi=\left( r_{11}\wedge \ldots \wedge r_{1s} \right) \vee \ldots \vee
\left( r_{l1} \wedge \ldots \wedge r_{ls}\right)
$$
and a variable $y_j$ belongs to a component $\Delta_t$. Then
$$
\eqalign{
\exists y_j \left( r_{11}\wedge \ldots \wedge r_{1s} \right) \vee
\ldots \vee \left( r_{l1} \wedge \ldots \wedge r_{ls}\right)
&=
\left( \exists y_j \left( r_{11}\wedge \ldots \wedge r_{1s} \right)\right)
\vee
\ldots \vee \left( \exists y_j\left( r_{l1} \wedge \ldots \wedge
r_{ls}\right)\right)
=\cr
&=
\left( r_{11}\wedge \ldots\wedge (\exists y_j r_{1 t})
\wedge\ldots\wedge r_{1s} \right) \vee
\ldots \vee \left( r_{l1} \wedge \ldots\wedge (\exists y_j) r_{l t}
\wedge \ldots \wedge r_{ls}\right),
}
$$
where the first equality is valid for any formulae and the second equality
is valid, because, in $i$th term of the disjunction, only $r_{it}$ may
depend on $y_j$. A similar transformation can be made for the
universal quantifier using the form $(C)$ of the formula $\psi\in\Lambda$.
Lemma is proven.

\medskip

Let us proceed with the proof of proposition.
By condition, each atomic subformula of~$\phi$
belongs to a class $\Phi_i$, i.e., in particular,
this subformula lies in $\Lambda$. Since
$\Lambda$ is closed with respect to logic operators and
quantifications on $y_j$ (by Lemma), we obtain that $\phi\in \Lambda$,
i.e.
$$
\phi = \bigvee_{i=1}^l\bigwedge_{j=1}^s r_{ij}.
$$
Let us assume that the maximum of $|\Delta_i|$ is attained at $\Delta_1$.
Then,
in each formula $r_{ij}$, where $j\neq 1$, we change the names
of bound variables from $\Delta_j$ for the names of variables from
$\Delta_1$. Since $\Delta_1$ is the largest component, we can
assign different names from $\Delta_1$ to different
variables (in each particular formula $r_{ij}$).
This renaming produces a formula $\phi'$ equivalent to $\phi$ and all
bound variables of $\phi'$ belong to $\Delta_1$.

%%%%%%%%%%%%%%%%%%%%%%%%%%%%%%%%%
\REFERENCES

\[Stru95]
Strunkov S. P.
On the theory of equations in finite groups
// Izv. Ross. Akad. Nauk. Ser. Mat. 1995.
V.59:6. P.171--180

\[Hall59]
%Холл М. Теория групп. - М.: ИЛ, 1962.
Hall M.
The theory of groups. New York, MacMillan Co., 1959.
13+434 pp.

\[AmV11]
Amit A., Vishne U.
Characters and solutions to equations in finite groups
// J. Pure Appl. Algebra. 2011. V.10. no.4. P.675--686.

\[AsTa01]
Asai T., Takegahara Y.
$|\Hom(A,G)|$, IV
//J. Algebra. 2001. 246. pp. 543--563.

\[BrTh88]
Brown K., Th\'evenaz J.
A generalization of Sylow's third theorem
//J. Algebra. 1988. 115. P. 414--430.

\[Frob03]
Frobenius G. \"Uber einen Fundamentalsatz der Gruppentheorie
// Berl. Sitz. 1903.  S.987--991.

\[GRV12]
Gordon C., Rodriguez-Villegas F.
On the divisibility of $\#\Hom(\Gamma, G)$ by $|G|$
// J. Algebra. 2012. V.350, no.1, P. 300--307.
See also arXiv:1105.6066.

\[Hall36]
Hall P.
On a theorem of Frobenius
// Proc. London Math. Soc. 1936. V.40. P.468--531.

\[Isaa70]
Isaacs I. M.
Systems of equations and generalized characters in groups
// Canad. J. Math. 1970. V.22. P.1040--1046.

\[Isaa12]
Isaacs I. M.
The number of group elements whose squares lie in a given subgroup
(an answer to Klyachko's question).
{\tt http://mathoverflow.net/questions/98639\#98809} (2012).

\[Iwa82]
S. Iwasaki,
A note on the $n$th roots ratio of a subgroup of a finite group
//
J. Algebra, 78:2 (1982), 460-474.
%ISSN 0021-8693, http://dx.doi.org/10.1016/0021-8693(82)90093-X.

\[Kula38]
Kulakoff A.
Einige Bemerkungen zur Arbeit: ``On a theorem of Frobenius" von P. Hall
// Matem. Sbornik. 1938. 3(45):2. 403--405.

\[Sehg62]
Sehgal S. K.
On P. Hall's generalisation of a theorem of Frobenius
//Proc. Glasgow Math. Assoc. 1962. 5.  P. 97--100.

\[Solo69]
Solomon L.
The solutions of equations in groups
// Arch. Math. 1969. V.20. no.3. P. 241--247.

\endREFERENCES
%%%%%%%%%%%%%%%%%%%%%%%%%%%%%%%%%
%%%%%%%%%%%%%%%%%%%%%%%%%%%%%%%%%
%%%%%%%%%%%%%%%%%%%%%%%%%%%%%%%%%

\end